\documentclass{article}

\input{diagrams}
\diagramstyle[PostScript=dvips,nohug]
\usepackage{epsfig}

\newtheorem{thm}{Theorem}[section]
\newtheorem{defn}[thm]{Definition}
\newtheorem{propn}[thm]{Proposition}

\newcommand{\mcm}[3]{\newcommand{#1}[#2]{{\ensuremath{#3}}}}
\usepackage{latexsym}
\mcm{\dashbk}{0}{\mbox{---}}
\mcm{\diagspace}{0}{\mbox{\hspace{2em}}}
\newcommand{\done}{\hfill\ensuremath{\Box}}

\newcommand{\bref}[1]{(\ref{#1})}
\newcommand{\ucontents}[2]{\addcontentsline{toc}{#1}{\numberline{}{#2}}}

\mcm{\cat}{1}{\mc{#1}}
\mcm{\fcat}{1}{\mb{#1}}
\mcm{\mc}{1}{\mathcal{#1}}
\mcm{\mr}{1}{\mathrm{#1}}
\mcm{\mb}{1}{\mathbf{#1}}
\newcommand{\url}[1]{\mbox{\tt #1 }}
\mcm{\eqv}{0}{\,\simeq\,}
\mcm{\iso}{0}{\,\cong\,}
\mcm{\of}{0}{\raisebox{0.2mm}{\ensuremath{\scriptstyle\circ}}}
\newcommand{\epsln}{\varepsilon}
\mcm{\op}{0}{\mr{op}}
\mcm{\Cat}{0}{\fcat{Cat}}
\mcm{\Mon}{0}{\mb{Mon}}
\mcm{\One}{0}{\fcat{1}}
\mcm{\Set}{0}{\fcat{Set}}
\mcm{\pr}{2}{\tuplebts{#1,#2}}
\mcm{\triple}{3}{\tuplebts{#1,#2,#3}}
\mcm{\homset}{3}{#1(#2,#3)}
\mcm{\go}{0}{\rTo}
\mcm{\goby}{1}{\rTo^{#1}}
\mcm{\goesto}{0}{\,\longmapsto\,}
\mcm{\goiso}{0}{\goby{\diso}}
\newarrow{Incl}C--->
\mcm{\diso}{0}{\sim}
\mcm{\Top}{0}{\fcat{Top}}
\mcm{\Topstar}{0}{\fcat{Top_*}}
\mcm{\ChCx}{0}{\fcat{ChCx}}
\mcm{\tuplebts}{1}{(#1)}
\mcm{\emm}{0}{\cat{M}}
\mcm{\Wzero}{0}{\scriptscriptstyle{\bullet}}
\mcm{\Wone}{0}{\raisebox{-.2em}{\epsfig{file=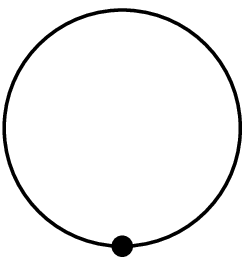,height=.9em}}}
\mcm{\Wtwo}{0}{\raisebox{-.2em}{\epsfig{file=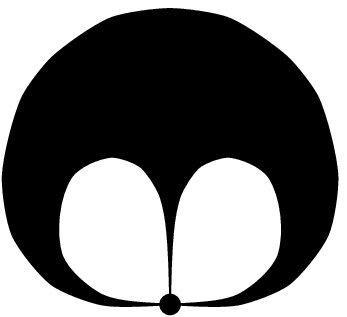,height=1.4em}}}
\mcm{\Wedgepic}{0}{\raisebox{-.2em}{\epsfig{file=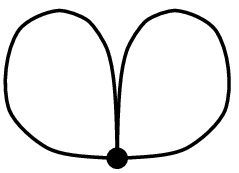,height=.7em}}}
\newenvironment{sketchpf}{\begin{trivlist}\item \textbf{Sketch
Proof}\ }{\end{trivlist}}
\mcm{\Loops}{0}{\homset{\Topstar}{S^1}{B}}

\begin{document}
\title{Up-to-Homotopy Monoids}
\author{Tom Leinster\\ \\
        \normalsize{Department of Pure Mathematics, University of
        Cambridge}\\ 
        \normalsize{Email: leinster@dpmms.cam.ac.uk}\\
        \normalsize{Web: http://www.dpmms.cam.ac.uk/$\sim$leinster}}
\date{}

\maketitle

\begin{abstract}

Informally, a homotopy monoid is a monoid-like structure in which properties
such as associativity only hold `up to homotopy' in some consistent way. This
short paper comprises a rigorous definition of homotopy monoid and a brief
analysis of some examples. It is a much-abbreviated version of the paper
`Homotopy Algebras for Operads' (\texttt{math.QA/0002180}), and
does not assume any knowledge of operads.

\end{abstract}

\tableofcontents

\section*{Introduction}
\ucontents{section}{Introduction}

This paper is a taster for the longer paper \cite{HAO}, `Homotopy algebras
for operads'. It is written for both the reader unfamiliar with the language
of operads, and the reader who would like to get the essence of \cite{HAO}
without taking too much time. The purpose here is to define what a homotopy
monoid is and to take a look at some examples.

I should first of all issue a warning that the term `homotopy monoid' is not
to be read in the same way as `homotopy group'. Informally, a homotopy monoid
is a monoid-like structure in which the associativity and unit laws don't
hold strictly, but only `weakly' or `up to homotopy' in some consistent
way. Two examples are well-known. In a (non-strict) monoidal category, one
has a tensor product $\otimes$ which is not (for instance) strictly
associative, but is associative up to isomorphism; moreover, these
isomorphisms obey rules of their own. The second example is that of a loop
space, that is, the space of all based loops in a fixed topological space
with basepoint. A loop space is not strictly speaking a topological monoid,
because of the fact that if one chooses a rule for composing any pair of
loops (e.g.\ travel each one at double speed) then that composition does not
obey unit or associative laws. However, these laws are obeyed up to
homotopy. Again, if one chooses particular homotopies to do this job then
these homotopies almost obey certain laws of their own---that is, they obey
them up to homotopy---and so on.

Those who do know about operads will know that monoids are the algebras for a
certain operad. In the full-length version \cite{HAO}, homotopy $P$-algebras
are defined for an arbitrary operad $P$, and the content of this paper is a
special case. Another special case is the operad $P$ whose algebras are
\emph{commutative} monoids; a homotopy $P$-algebra in the category of
topological spaces is then exactly what Segal defined as a $\Gamma$-space in
\cite{Seg} (or `special $\Gamma$-space', in more recent
terminology). It is from Segal's paper that this work is descended.

\subsection*{Acknowledgements}

I am very grateful to the many people who have helped me with this research;
a full-length acknowledgement is contained in the full-length version
\cite{HAO} of this paper. I would like to add my thanks here to two further
parties: firstly, to the organisers of the 71st Peripatetic Seminar on
Sheaves and Logic at Louvain-la-Neuve, 16--17 October 1999, who gave me the
opportunity to present this work and therefore the stimulus to write it up
in this form, and secondly, to George Janelidze, who alerted me to the
ambiguity of the phrase `homotopy monoid'.

The work was supported by the Laurence Goddard Fellowship at St John's
College, Cambridge. The arrows in the document were generated using Paul
Taylor's diagrams package.

\section{The Environment}	\label{sec:env}

Let \emm\ be a monoidal category. Our aim is to define what a homotopy monoid
in \emm\ is, and so we need some notion of `homotopy' or `weak
equivalence' in \emm. This will be achieved by a very simple device: we
mark out certain of the morphisms in \emm\ and decide to call them
homotopy equivalences, or just `equivalences'.

So: from now on suppose that \emm\ is equipped with a distinguished
subcollection of its morphisms, and call these morphisms the \emph{(homotopy)
equivalences}\/ in \emm. For reasons which are more apparent in \cite{HAO},
we will only consider examples having the following properties:
\begin{itemize}
\item any isomorphism is an equivalence
\item if $h=g\of f$ is a composite of morphisms in \emm, and if any two of
$f$, $g$, $h$ are equivalences, then so is the third
\item if $A\goby{f}B$ and $A'\goby{f'}B'$ are equivalences then so is
$A\otimes A' \goby{f\otimes f'}$ $B \otimes B'$.
\end{itemize}

\subsubsection*{Examples}

\begin{enumerate}
\item
\emm\ is the monoidal category \pr{\Top}{\times} of topological spaces with
the usual (cartesian) product; homotopy equivalences are homotopy
equivalences.
\item
\emm\ is the monoidal category \pr{\Cat}{\times} of (small) categories and
functors, with cartesian product; equivalences are equivalences of
categories---that is, those functors $G$ for which there exists some functor
$F$ with $G\of F \iso 1$ and $F\of G \iso 1$. ($F$ will be called a
\emph{pseudo-inverse} to $G$.)
\item	\label{eg:isos}
\emm\ is any monoidal category, and the equivalences are just the
isomorphisms. Given the first of the three conditions above, this is the
smallest possible class of equivalences in \emm.
\item
\emm\ is \pr{\ChCx}{\otimes}, the monoidal category of chain complexes (of
modules over a fixed commutative ring) with the usual tensor product;
equivalences are chain homotopy equivalences.
\item
\emm\ is \pr{\Topstar}{\vee}, where \Topstar\ is the category of spaces
with basepoint and continuous basepoint-preserving maps, and $\vee$ is the
wedge product (join two spaces together by their basepoints); equivalences
are homotopy equivalences which respect basepoints.
\end{enumerate}

\section{The Definition of Homotopy Monoid}	\label{sec:defn}

We are nearly ready to give the definition of homotopy monoid, but first need
a preliminary definition; it is one of the various notions of a map between
monoidal categories.

\begin{defn}
Let \cat{N'} and \cat{N} be monoidal categories. A \emph{colax monoidal
functor} $\cat{N'} \go \cat{N}$ consists of a functor $X: \cat{N'} \go
\cat{N}$ together with maps 
\[
\begin{array}{rccc}
\xi_{AB}:	&X(A\otimes B) 	&\go  	&X(A) \otimes X(B)	\\
\xi_0:		&X(I) 		&\go 	&I
\end{array}
\]
($A,B \in \cat{N'}$) satisfying naturality and coherence axioms.
\end{defn}

Here $\otimes$ and $I$ denote the tensor operation and unit object in both
monoidal categories \cat{N'} and \cat{N}. The coherence axioms can be found
in \cite[1.1.1]{HAO} or in \cite[6.4.1]{Borc}. If using the latter source one
must reverse all the arrows labelled by a $\tau$, since what is being defined
there is a \emph{lax}\/ monoidal functor.

\begin{defn}
Let \emm\ be a monoidal category equipped with a class of equivalences, as in
Section~\ref{sec:env}. A \emph{homotopy monoid}\/ in \emm\ is a colax
monoidal functor
\[
\pr{X}{\xi}: \Delta \go \emm
\]
for which the maps $\xi_0$ and $\xi_{mn}$ are each homotopy equivalences
($m,n\geq 0$).
\end{defn}

Here $\Delta$ is the category of finite ordinals, equivalent to the category
of totally ordered finite sets (including $\emptyset$). The monoidal product
on $\Delta$ is addition, and the unit object is the empty set. The objects of
$\Delta$ are written $0$, $1$, \ldots, so that $n$ denotes an $n$-element
totally ordered set.

The rest of this paper is devoted to examples of homotopy monoids. The
non-trivial examples take a while to explain, so this section only covers
them in a very sketchy way; we come back to two of them in
Sections~\ref{sec:moncats} and~\ref{sec:loops}. Further details on all can be
found in \cite{HAO}.

\subsubsection*{Examples}

\begin{enumerate}

\item \label{eg:iso-mon} 
Suppose that the equivalences in \emm\ are just the isomorphisms, as in
Example \ref{sec:env}\bref{eg:isos}. Then a homotopy monoid in \emm\ is just
a monoidal functor $\Delta\go\emm$. By a monoidal functor I mean what is
sometimes called a strong monoidal functor, that is, a map of monoidal
categories which preserves the tensor and unit up to coherent
isomorphism. But it is well-known that monoidal functors $\Delta\go\emm$ are
essentially just monoids in \emm. Thus when \emm\ has no `interesting'
homotopy equivalences, a homotopy monoid in \emm\ is merely a monoid in \emm.

To be a little more precise, it is proved (essentially) in
\cite[VII.5.1]{CWM} that the category of monoidal functors and monoidal
transformations from $\Delta$ to \emm\ is equivalent to the category of
monoids in \emm. A monoidal functor \pr{X}{\xi} corresponds to a certain
monoid structure on the object $X(1)$ of \emm, and in general we will regard
$X(1)$ as being in some sense the underlying space of a homotopy monoid
\pr{X}{\xi}.

\item
A genuine monoid in \Cat\ is a strict monoidal category; one would therefore
expect a homotopy monoid in \Cat\ to be something comparable to a
(non-strict) monoidal category. This is discussed further in
Section~\ref{sec:moncats}.

\item
Any loop space is a homotopy monoid in \Top: see Section~\ref{sec:loops}.

\item
A monoid in \ChCx\ is a differential graded algebra, so a homotopy monoid in
\ChCx\ might be called a homotopy differential graded algebra. But there is
already a commonly-used formulation of the concept of homotopy differential
graded algebra, namely $A_\infty$-algebras. So homotopy monoids in \ChCx\
ought to be related somehow to $A_\infty$-algebras; this is not discussed
further here, but is in \cite[3.5]{HAO}. Similar remarks apply to
$A_\infty$-spaces and homotopy monoids in \Top.

\end{enumerate}

\section{Monoidal Categories}	\label{sec:moncats}

Let us now look more closely at homotopy monoids in \Cat. Such a structure
consists of a functor $C:\Delta \go \Cat$ (previously called $X$) together
with equivalences of categories
\[
\begin{array}{rccc}
\xi_{mn}:	&C(m+n) 	&\go 	&C(m) \times C(n)	\\
\xi_{0}:	&C(0) 		&\go 	&\One
\end{array}
\]
($m,n\geq 0$) fitting together nicely. Note that by assembling the components
of $\xi$ we obtain a (canonical) equivalence $C(n) \goby{\eqv} C(1)^n$ for
each $n$; we regard $C(1)$ as the `base category' of the homotopy monoid
\pr{C}{\xi}.
% , in a sense made more clear by the 
% proof of Proposition~\ref{propn:hty-mon-cat}.  

There now follows a comparison (in one direction) between homotopy monoids in
\Cat\ and monoidal categories in the usual sense.

\begin{propn}	\label{propn:hty-mon-cat}
A homotopy monoid in \Cat\ gives rise to a monoidal category.
\end{propn}

\begin{sketchpf}
Take a homotopy monoid \pr{C}{\xi} in \Cat, and construct from it a monoidal
category as follows.

\begin{description}
\item[Underlying category:] $C(1)$.

\item[Tensor:] What we want to define is a functor
\[
\otimes: C(1) \times C(1) \go C(1);
\]
what we actually have are functors
\begin{diagram}[size=2.5em]
		&			&C(2)	&		&	\\
		&\ldTo<{\xi_{1,1}}>{\eqv}&	&\rdTo>{C(!)}	&	\\
C(1)\times C(1)	&			&	&		&C(1)	\\
\end{diagram}
where $!$ is the unique map $2\go 1$ in $\Delta$. So for each $m$ and $n$,
choose (arbitrarily) a pseudo-inverse $\psi_{mn}$ to $\xi_{mn}$, and define
$\otimes$ as the composite
\[
C(1) \times C(1) \goby{\psi_{1,1}} C(2) \goby{C(!)} C(1).
\]

\item[Associativity isomorphisms:] The next piece of data we need is a
natural isomorphism between $\otimes\of(\otimes\times 1)$ and
$\otimes\of(1\times\otimes)$. To see why such an isomorphism should exist,
consider what would happen if the $\psi_{mn}$'s were \emph{genuine} inverses
to the $\xi_{mn}$'s. Then the $\psi_{mn}$'s would satisfy the same coherence
and naturality axioms as the $\xi_{mn}$'s (with the arrows reversed), and
this would guarantee that all sensible diagrams built up out of $\psi_{mn}$'s
commuted. Hence $\otimes$ would be strictly associative. As it is,
$\psi_{mn}$ is only inverse to $\xi_{mn}$ up to isomorphism, and
correspondingly $\otimes$ is associative up to isomorphism.

In practice, choose (at random) natural isomorphisms
\[
\eta_{mn}: 1 \goiso \xi_{mn}\of\psi_{mn},
\diagspace
\epsln_{mn}: \psi_{mn}\of\xi_{mn} \goiso 1
\]
for each $m$ and $n$. Then a natural isomorphism
\[
\alpha: \otimes\of(\otimes\times 1) \goiso \otimes\of(1\times\otimes)
\]
can be built up from the $\eta_{mn}$'s and $\epsln_{mn}$'s.

\item[Pentagon:] We must now check that the associativity isomorphism just
defined satisfies the famous pentagon coherence axiom. This asserts the
commutativity of a certain diagram built up from components of $\alpha$, that
is, built up from $\eta_{mn}$'s and $\epsln_{mn}$'s. However, this diagram
does \emph{not} commute, which is perhaps unsurprising since $\eta_{mn}$ and
$\epsln_{mn}$ were chosen independently.

But all is not lost: for recall the result that if
\[
\begin{array}{ll}
F: A \go B,			&G: B\go A,		\\
\sigma: 1 \goiso G\of F,	&\tau: F\of G \goiso 1
\end{array}
\]
is an equivalence of categories, then $\tau$ can be modified to another
natural isomorphism $\tau'$ so that $(F,G,\sigma,\tau')$ is both an
adjunction and an equivalence (see \cite[IV.4.1]{CWM}). So when we chose the
natural isomorphisms $\eta_{mn}$ and $\epsln_{mn}$ above, we could have done
it so that $(\psi_{mn},\xi_{mn},\eta_{mn},\epsln_{mn})$ was an
adjunction. Assume that we did so. Then this being an adjunction says that
certain basic diagrams involving $\eta_{mn}$ and $\epsln_{mn}$ commute
(namely, the diagrams for the triangle identities \cite[IV.1]{CWM}): and that
is enough to ensure that the pentagon commutes.

\item[Units:] We also need to define the unit object for the tensor and the
left and right unit isomorphisms, and to prove those coherence conditions
which involve units. This is done by the same methods as above.
\done
\end{description}
\end{sketchpf}

So starting from a homotopy monoid in \Cat, we have constructed a monoidal
category. The construction involves arbitrary choices, but these choices only
affect the resulting monoidal category up to isomorphism (that is,
isomorphism in the category of monoidal categories and monoidal functors).

\section{Loop Spaces}	\label{sec:loops}

So far we have not actually seen any examples of homotopy monoids beyond the
trivial. We remedy this now by sketching an argument that any loop space is a
homotopy monoid in \Top. Liberal use will be made of function spaces
(exponentials); the conscientious reader should therefore regard `space' as
meaning `compactly generated Hausdorff space', or take some other convenient
cartesian closed substitute for the category of topological spaces.

Fix a space $B$ with basepoint and form the loop space of $B$, that is, the
space \Loops\ of continuous basepoint-preserving maps from the circle $S^1$
to $B$.

\begin{propn}
\Loops\ is a homotopy monoid in \Top, in a canonical way.
\end{propn}

\begin{sketchpf}
The object is to find a homotopy monoid \pr{X}{\xi} in \Top\ with
$X(1)=\Loops$. (To see why this is a reasonable interpretation of the
statement of the Proposition, compare the role of $X(1)$ in
Example~\ref{sec:defn}\bref{eg:iso-mon} and $C(1)$ in
Section~\ref{sec:moncats}.)  The proof is in two steps.

\begin{enumerate}
\item	\label{part:comon}
$S^1$ is a homotopy comonoid, in the sense that there is a colax monoidal
functor 
\[
\pr{W}{\omega}: \pr{\Delta}{+} \go \pr{\Top_{*}^\op}{\vee}
\]
with $W(1)=S^1$ and with the components $\omega_0$, $\omega_{mn}$ of $\omega$
all homotopy equivalences. In other words, $S^1$ is a homotopy monoid in
\pr{\Top_{*}^\op}{\vee}. We define \pr{W}{\omega} as follows:
\begin{itemize}
\item
$W(n)$ is the standard $n$-simplex $\Delta^n$ with its $(n+1)$ vertices
collapsed to a single point, and this declared the basepoint, e.g.\ 
\begin{eqnarray*}
W(0)	&=	&\Wzero\ ,	\\
W(1)	&=	&\Wone = S^1,	\\
W(2)	&=	&\Wtwo.
\end{eqnarray*}
\item
$W$ is defined on morphisms using the standard face and degeneracy maps of
simplices
\item
$\omega$ is defined by face maps: for instance,
\[
\omega_{1,1}: W(1) \vee W(1) \go W(2)
\]
is the inclusion
\[
\Wedgepic \rIncl \Wtwo,
\]
which is evidently a homotopy equivalence. ($W(2)$ can be thought of as a
thickened wedge of two circles.)
\end{itemize}

\item	\label{part:repble}
There's a functor
\[
\homset{\Topstar}{\dashbk}{B}: 
\pr{\Top_{*}^\op}{\vee} \go \pr{\Top}{\times},
\]
which is a (strong) monoidal functor (in the terminology of
Example~\ref{sec:defn}\bref{eg:iso-mon}) simply because $\vee$ is the
coproduct in $\Topstar$. Note also that this functor preserves homotopy
equivalences.
\end{enumerate}

Composing the functors of \bref{part:comon} and \bref{part:repble} yields a
homotopy monoid $\pr{X}{\xi}: \pr{\Delta}{+} \go \pr{\Top}{\times}$ with
$X(1) = \Loops$, as required.
\done
\end{sketchpf}

Let us finish by examining the homotopy monoid structure we have just put on
the loop space, and in particular at how the composition of two loops is
handled. 

We have 
\[
X(2) = \homset{\Topstar}{\Wtwo}{B},
\]
and the pieces of \pr{X}{\xi} relevant to binary composition are the maps
\begin{diagram}[size=2.5em]
	&				&X(2)	&		&	\\
	&\ldTo<{\xi_{1,1}}>{\eqv}	&	&\rdTo>{X(!)}	&	\\
X(1)^2	&				&	&		&X(1).	\\
\end{diagram}
This diagram is
\begin{equation}	\label{eqn:canonical}
\begin{diagram}[size=2.5em]
	&		&\homset{\Topstar}{\Wtwo}{B}	&	&	\\
	&\ldTo<{\eqv}	&				&\rdTo	&	\\
\homset{\Topstar}{\Wedgepic}{B}	= \Loops^2&	&	&	&\Loops,\\
\end{diagram}
\end{equation}
where the map on the left is restriction to the two inner circles and the map
on the right is restriction to the outer circle. Note that all of the data
making up \pr{X}{\xi}, and in particular the maps in~\bref{eqn:canonical},
is constructed canonically from $B$: no arbitrary choices have been made. In
contrast, there is no canonical map
\[
\Loops^2 \go \Loops
\]
defining `composition': although the obvious and customary choice is to use
the map described by the instruction `travel each loop at double speed', this
appears to have no particular advantage or special algebraic status compared
to any other choice. Since the usual formulation of the idea of homotopy
topological monoid, $A_\infty$-spaces, does entail this arbitrary choice of a
composition law for loops, one might regard this as a virtue of the
definition presented here.

\end{document}